\documentclass{amsart}
\usepackage{wasowicz,graphicx}

\newtheorem*{thA}{Theorem A}
\newtheorem*{thB}{Theorem B}
\newtheorem*{thC}{Theorem C}
\newtheorem*{thD}{Theorem D}

\info{%
 Published: J. Math. Anal. Appl. 332~(2007), 1229--1241\\
 In this version of the manuscript the bad example of Remark~\ref{bad:ex} is corrected.}

\begin{document}

\title[Generalized convexity: support properties, Hadamard inequalities]%
      {Support--type properties of convex functions of higher order and Hadamard--type inequalities}
\author{\SW}
\address{\SWaddr}
\email{\SWmail}
\date{July 25, 2008}
\keywords{%
 Approximate integration,
 error bounds,
 Hadamard inequality,
 higher--order convexity,
 quadrature rules,
 support theorems}
\subjclass[2000]{Primary: 26A51, Secondary: 26D15, 41A55, 41A80}

\begin{abstract}
 It is well--known that every convex function $f:I\to\R$ (where
 $I\subset\R$ is an interval) admits an affine support at every
 interior point of~$I$ (i.e. for any $x_0\in\Int I$ there exists an
 affine function $a:I\to\R$ such that $a(x_0)=f(x_0)$ and $a\le f$
 on~$I$). Convex functions of higher order (precisely of an odd
 order) have a~similar property: they are supported by the
 polynomials of degree no greater than the order of convexity. In this
 paper the attaching method is developed. It is applied to
 obtain the general result --- Theorem~\ref{attaching}, from which the
 mentioned above support theorem and some related properties of convex
 functions of higher (both odd and even) order are derived. They are
 applied to obtain some known and new Hadamard--type inequalities
 between the quadrature operators and the integral approximated by
 them. It is also shown that the error bounds of quadrature rules
 follow by inequalities of this kind.
\end{abstract}
\maketitle

\section{Introduction}

Let $I\subset\R$ be an interval and $f:I\to\R$. For distinct
points of~$I$ the \emph{divided differences} of~$f$ are defined
recursively as follows: $[x_1,f]:=f(x_1)$ and
\begin{equation}\label{eq:div_diff_def}
 [x_1,\dots,x_{n+1};f]:=\frac{[x_2,\dots,x_{n+1};f]
                       -[x_1,\dots,x_n;f]}{x_{n+1}-x_1},\quad
 n\in\N,\;n\ge 2.
\end{equation}
For $n$ distinct points $x_1,\dots,x_n\in I$ ($n\ge 2$) the following
formula holds true:
\begin{equation}\label{eq:div_diff}
 [x_1,\dots,x_n;f]=\frac{D(x_1,\dots,x_n;f)}{V(x_1,\dots,x_n)},
\end{equation}
where
\[
 D(x_1,\dots,x_n;f):=
 \begin{vmatrix}
  1&\dots&1\\
  x_1&\dots&x_n\\
  \vdots&&\vdots\\
  x_1^{n-2}&\dots&x_n^{n-2}\\
  f(x_1)&\dots&f(x_n)
 \end{vmatrix}
\]
and $V(x_1,\dots,x_n)$ stands for the Vandermonde determinant of the
terms involved. By~\eqref{eq:div_diff} we can immediately see that the
divided differences are symmetric.
\par
Let $\Pi_n$ be the family of all polynomials of degree at most~$n$.
For $k$ distinct points $x_1,\dots,x_k\in I$ denote by
$P(x_1,\dots,x_k;f)$ the (unique) interpolation polynomial
$p\in\Pi_{k-1}$ such that $p(x_i)=f(x_i)$, $i=1,\dots,k$. Then for any
$n+1$ distinct points $x_1,\dots,$ $x_{n+1}\in I$ and for any
$x\in I\setminus\{x_1,\dots,x_{n+1}\}$ we have
\begin{equation}\label{eq:poly}
 f(x)-P(x_1,\dots,x_{n+1};f)(x)=[x_1,\dots,x_{n+1},x;f]
                                \prod_{i=1}^{n+1}(x-x_i).
\end{equation}
For the definition and properties of divided differences the reader is
referred to~\cite{K,P,RV}.
\par
If $x_1,\dots,x_{n+1}\in I$ are distinct then Newton's
Interpolation Formula holds:
\begin{multline}\label{eq:Newton}
 P(x_1,\dots,x_{n+1};f)(x)=f(x_1)+[x_1,x_2;f](x-x_1)+\dots\\
    +[x_1,\dots,x_{n+1};f](x-x_1)\dotsm(x-x_n).
\end{multline}
Next we recall the notion of convex functions of higher order.
Hopf's thesis~\cite{H} from 1926 seems to be the first work devoted to
this topic (the functions with nonnegative divided differences
were considered but the name ``convex functions of higher order''
was not used). Eight years later higher--order convexity was
extensively studied by Popoviciu~\cite{P} (cf. also~\cite{K,RV}).
Let $n\in\N$. A function $f:I\to\R$ is called \emph{$n$--convex\/}
if $[x_1,\dots,x_{n+2};f]\ge 0$ for any $n+2$ distinct points
$x_1,\dots,x_{n+2}\in I$. It follows by~\eqref{eq:div_diff} that $f$
is $n$--convex if and only if
\begin{equation}\label{ineq:nconv}
 D(x_1,\dots,x_{n+2};f)\ge 0
\end{equation}
for any $x_1,\dots,x_{n+2}\in I$ with $x_1<\dots<x_{n+2}$ (since
$V(x_1,\dots,x_{n+2})>0$).
\par
For $n=1$ it is not difficult to observe that the $n$--convexity
reduces to convexity in the usual sense.
\par
By~\eqref{eq:poly} we obtain the following important property of
convex functions of higher order (cf.~\cite{K,P,RV}): a function
$f:I\to\R$ is $n$--convex if and only if for any
$x_1,\dots,x_{n+1}\in I$ with $x_1<\dots<x_{n+1}$ the graph
of an interpolation polynomial $p:=P(x_1,\dots,x_{n+1};f)$ passing
through the points $\bigl(x_i,f(x_i)\bigr)$, $i=1,\dots,n+1$, changes
succesively the side of the graph of~$f$ (always $p(x)\le f(x)$ for
$x\in I$ such that $x>x_{n+1}$, if such points do exist). More
precisely,
\begin{align}\label{ineq:nconv_poly}
 \nonumber
 (-1)^{n+1}\bigl(f(x)-p(x)\bigr)&\ge 0,\quad x<x_1,\;x\in I,\\
 (-1)^{n+1-i}\bigl(f(x)-p(x)\bigr)&\ge 0,
 \quad x_i<x<x_{i+1},\;i=1,\dots n,\\
 \nonumber
 f(x)-p(x)&\ge 0,\quad x>x_{n+1},\;x\in I.
\end{align}
The theorem below contains another property of higher--order convexity
(cf.~\cite[p.~391, Corollary~1]{K}, \cite[p.~27]{P}).
\begin{thA}
 If $f:I\to\R$ is $n$--convex then for any $k\in\{1,\dots,$ $n+1\}$
 the divided differences $[x_1,\dots,x_k;f]$ are bounded on every
 compact interval $[a,b]\subset\Int I$.
\end{thA}
Convex functions of higher order have the following regularity
property (cf. \cite{BP_Cheb,K,P}):
\begin{thB}
 If $f:[a,b]\to\R$ is $n$--convex then $f$ is continuous on $(a,b)$
 and bounded on $[a,b]$.
\end{thB}
All the integrals that appear in this paper are understood in the
sense of Riemann. Then by Theorem~B we obtain
\begin{thC}
 If $f:[a,b]\to\R$ is $n$--convex then $f$ is integrable on $[a,b]$.
\end{thC}
For $n$--convex functions which are ($n+1$)--times differentiable the
following result holds (cf.~\cite[Theorems A and B]{W1},
\cite[Theorems 1.2 and 1.3]{W2}, cf. also \cite{K,P,RV}):
\begin{thD}
 Assume that $f:[a,b]\to\R$ is ($n+1$)--times differentiable on
 $(a,b)$ and continuous on $[a,b]$. Then~$f$ is $n$--convex if and
 only if $f^{(n+1)}(x)\ge 0$, $x\in(a,b)$.
\end{thD}
It is well--known that every convex function $f:I\to\R$ admits an
affine support at every interior point of~$I$ (i.e. for any
$x_0\in\Int I$ there exists an affine function $a:I\to\R$ such that
$a(x_0)=f(x_0)$ and $a\le f$ on~$I$). Convex functions of higher
order (precisely of an odd order) have a similar property: they are
supported by the polynomials of degree no greater than the order of
convexity. Such a~result was obtained by Ger~\cite{G}, who
assumed that the supported $n$--convex function, defined on an open
and convex subset of a~normed space, was of the class $\C^{n+1}$. In
this paper we develop the attaching method and we use it to
prove in Theorem~\ref{attaching} a~support--type result of a general
nature. As almost immediate consequences we obtain the result
improving Ger's theorem (we remove the differentiability assumption)
for functions defined on a real interval and more support--type
properties of convex functions of higher (both odd and even) order.
\par
In the theory of convex functions an important role is played by the
famous Hermite--Hadamard inequality. It states that if $f:[a,b]\to\R$
is convex then
\begin{equation}\label{ineq:Hadamard}
 f\Bigl(\frac{a+b}{2}\Bigr)\le
 \frac{1}{b-a}\int_a^b f(x)dx\le\frac{f(a)+f(b)}{2}.
\end{equation}
The interesting study of this inequality and lots of related
inequalities was given by Dra\-go\-mir and Pearce~\cite{DP}. In this paper
we apply the above mentioned support--type properties of convex
functions of higher order to obtain both known and new Hadamard--type
inequalities between the quadrature operators and the integral
approximated by them. We also show that the error bounds of quadrature
rules follow by inequalities of this kind.

\section{Attaching method}

\begin{figure}[h]
 \begin{center}
  \includegraphics{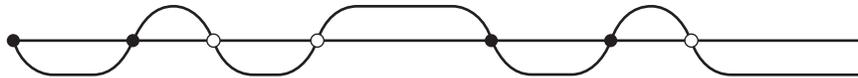}
  \caption{The ``bubbles'' are attached to the ``bullets''.}
  \label{Fig1}
 \end{center}
\end{figure}

In this section we describe this method. Consider the $n$--convex
function $f:I\to\R$ and take the polynomial $p\in\Pi_n$
interpolating~$f$ at $n+1$ distinct points of~$I$. The
Figure~\ref{Fig1} is drawn for the $6$--convex function. The graph
of~$f$ is represented by the horizontal straight line. Then by
$n$--convexity the graph of $p$ (symbolized by the curve line) meeting
the graph of~$f$ changes successively its side.

\par
\begin{figure}[h]
 \begin{center}
  \includegraphics{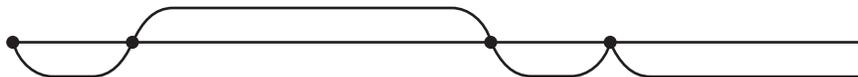}
  \caption{The degree (its upper bound) of the polynomial is preserved.}
  \label{Fig2}
 \end{center}
\end{figure}

If the ``bubbles'' lying between two consecutive ``bullets''
tend to the nearest ``bullet'' situated on the left hand
side of them, then the appropriate interpolation polynomials belonging
to $\Pi_n$ (the graph of one of them is shown at the
Figure~\ref{Fig1}) tend to some polynomial belonging to $\Pi_n$. Then
we arrive at the situation that we can see at the Figure~\ref{Fig2}.
\par
Figures~\ref{Fig1} and~\ref{Fig2} have only an explanatory
character. The curve lines are not really the graphs of polynomials
and the straight line is not really the graph of~$f$. They
illustrate only the location of graphs of appropriate
polynomials on a proper side of the graph of $f$.

\section{Support--type theorem}

Now we are going to prove a support--type theorem of the general
nature. In the proof we would like to use the boundedness of
divided differences of an $n$--convex function. By Theorem~A this is
the case when all the points involved belong to the compact
subinterval of $\Int I$. However, we need also this property for
divided differences involving  additionally the boundary points of $I$
(if they do exist). That is why we prove below the following lemma.
\begin{lem}\label{bounded_cl}
 Let $n\in\N$, $A\subset\R$, $a\not\in\Cl A$ and $f:A\cup\{a\}\to\R$.
 If for any $k\in\{1,\dots,n\}$ the divided differences
 $[x_1,\dots,x_k;f]$ are bounded on~$A$ then for any
 $k\in\{1,\dots,n\}$ they remain bounded on $A\cup\{a\}$.
\end{lem}
\begin{proof}
 For $k=1$ there is nothing to prove. For $k>1$ assume that the
 divided differences $[x_1,\dots,x_{k-1};f]$ are bounded on
 $A\cup\{a\}$. To finish the proof it is enough to show the assertion
 for $k$--point divided differences containing~$a$. To proceed this
 job take $x_1,\dots,x_{k-1}\in A$. Then by~\eqref{eq:div_diff_def}
 \[
  \bigl|[a,x_1,\dots,x_{k-1};f]\bigr|=
  \frac{\bigl|[x_1,\dots,x_{k-1};f]-[a,x_1,\dots,x_{k-2};f]\bigr|}
       {|x_{k-1}-a|}\le\frac{2M}{|a-A|},
 \]
 where
 \begin{align*}
  M&:=\sup\Bigl\{\bigl|[x_1,\dots,x_{k-1};f]\bigr|\;:\;
      x_1,\dots,x_{k-1}\in A\cup\{a\}\Bigr\},\\
  |a-A|&:=\inf\bigl\{\,|a-b|\;:\;b\in A\,\bigr\}>0
  \quad\text{(since $a\not\in\Cl A$).}
 \end{align*}
\end{proof}
Now we are ready to prove the main result of this section.
\begin{thm}\label{attaching}
 Let $n\in\N$ and $f:I\to\R$ be an $n$--convex function. Fix $k\in\N$,
 $k\le n$ and take $x_1,\dots,x_k\in I$ such that $x_1<\dots<x_k$.
 Assign to each point $x_j$ ($j=1,\dots,k$) the multiplicity
 $l_j\in\N$ ($l_j-1$ stands for the number of points attached to
 $x_j$). We require $l_1+\dots+l_k=n+1$ and if $x_1=\inf I$ then
 $l_1=1$, if $x_k=\sup I$ then $l_k=1$ (the points can be attached
 only to the interior points of~$I$). Denote $I_0=(-\infty,x_1)$,
 $I_j=(x_j,x_{j+1})$, $j=1,\dots,k-1$ and $I_k=(x_k,\infty)$. Under
 these assumptions there exists a~polynomial $p\in\Pi_n$ such that
 $p(x_j)=f(x_j)$, $j=1,\dots,k$ and
 \begin{align}\label{ineq:attaching}
  \nonumber
  (-1)^{n+1}\bigl(f(x)-p(x)\bigr)&
  \ge 0\quad\text{for }x\in I_0\cap I,\\
  (-1)^{n+1-(l_1+\dots+l_j)}\bigl(f(x)-p(x)\bigr)&
  \ge 0\quad\text{for }x\in I_j,\;j=1,\dots,k-1,\\
  \nonumber
  f(x)-p(x)&\ge 0\quad\text{for }x\in I_k\cap I.
 \end{align}
\end{thm}
Before we start the proof let us notice that at the Figures~\ref{Fig1}
and~\ref{Fig2} we have $n=6$, $k=4$ (the ``bubbles'' were attached to
the ``bullets'') and the multiplicities of the ``bullets''
are (from left to right) 1, 3, 1, 2, respectively.
\begin{proof}[Proof of Theorem~\ref{attaching}.]
 Let $m\in\N$. If $l_j>1$ ($j=1,\dots,k$), we take the points
  \[
   x_j<x_j+\frac{1}{m}<\dots<x_j+\frac{l_j-1}{m}.
  \]
  For $m$ large enough all these points belong to $I_j$. Then the
  sequence
  \begin{equation}\label{*}
   \Bigl(
    x_1,x_1+\frac{1}{m},\dots,x_1+\frac{l_1-1}{m},x_2,
    \dots,x_k,x_k+\frac{1}{m},\dots,x_k+\frac{l_k-1}{m}
   \Bigr)
  \end{equation}
  is increasingly ordered and it contains $n+1$ points of $I$ (because
  of $l_1+\dots+l_k=n+1$). There exists a~polynomial $p_m\in\Pi_n$
  interpolating~$f$ at the points of the sequence~\eqref{*}. We use
  Newton's Interpolation Formula~\eqref{eq:Newton} to write $p_m$.
  This formula contains the products of binomials of the form
  \[
   x-x_j-\frac{s_j}{m},\qquad j=1,\dots,k,
   \quad s_j=0,\dots,l_j-1\;(k\le n)
  \]
  and the divided differences involving points of the
  sequence~\eqref{*}. By Theorem~A and Lemma~\ref{bounded_cl} the
  sequences of these differences containing 1~point, 2~points, \dots,
  $n+1$ points, respectively, are bounded and for that reason they
  contain the convergent subsequences. By taking (if needed) the
  common subsequence $(\alpha_m)$ of positive integers we may assume
  without loss of generality that all these sequences are convergent.
  \par
  Let $x\in I$ and $p(x):=\lim\limits_{m\to\infty} p_m(x)$.
  Then $p\in\Pi_n$ and by the construction we have
  $p(x_j)=f(x_j)$, $j=1,\dots,k$.
  \par
  Let $x\in I_0\cap I$ (if $x\in I_0\cap I\ne\varnothing$). Then
  $x<x_1$ and by $n$--convexity and~\eqref{ineq:nconv_poly}
  $(-1)^{n+1}\bigl(f(x)-p_m(x)\bigr)\ge 0$ for $m$ large enough.
  Tending with~$m$ to infinity we get
  $(-1)^{n+1}\bigl(f(x)-p(x)\bigr)\ge 0$.
  \par
  Let $x\in I_k\cap I$ (if $x\in I_k\cap I\ne\varnothing$). Then for
  $m$ large enough we have
  \[
   x>x_k+\frac{l_k-1}{m}.
  \]
  We infer by~\eqref{ineq:nconv_poly} that $f(x)-p_m(x)\ge 0$, whence
  letting $m\to\infty$ we obtain $f(x)-p(x)\ge 0$.
  \par
  Finally let $x\in I_j$, $j=1,\dots,k-1$. For $m$ large enough we
  have
  \[
   x_j+\frac{l_j-1}{m}<x<x_{j+1}.
  \]
  Observe that the point $x_j+\frac{l_j-1}{m}$ has in the
  sequence~\eqref{*} the number $l_1+\dots+l_j$. Therefore
  by~\eqref{ineq:nconv_poly}
  $(-1)^{n+1-(l_1+\dots+l_j)}\bigl(f(x)-p_m(x)\bigr)\ge 0$. For
  $m\to\infty$ we get
  \[
   (-1)^{n+1-(l_1+\dots+l_j)}\bigl(f(x)-p(x)\bigr)\ge 0,
  \]
  which finishes the proof.
\end{proof}
\begin{rem}
  In the classical setting, if $f:I\to\R$ admits at each point
  $x_0\in\Int I$ an affine support, then $f$ is convex. This is also
  the case for the statement of Theorem~\ref{attaching}: it
  characterizes $n$--convexity. Indeed, to prove that $f$ is
  $n$--convex it is enough to assume that the appropriate polynomial
  exists for $k=n$, $x_1,\dots,x_n\in\Int I$ with $x_1<\dots<x_n$
  and $l_1=\dots=l_{n-1}=1$, $l_n=2$. This is shown by the present
  author in \cite[Theorem 3]{Was_Cheb} in a~more general setting, i.e.
  for convex functions with respect to Chebyshev systems (for
  a~polynomial Chebyshev system $(1,x,\dots,x^n)$ such a~convexity
  reduces to $n$--convexity).
  We have formulated Theorem~\ref{attaching} in the form of the
  necessary condition since, as we can see, the sufficient condition
  can be weakened.
\end{rem}
\begin{rem}\label{bad:ex}
  The polynomial obtained in Theorem~\ref{attaching} need not to be
  unique. Let $n=5$, $k=3$. Then for a~5--convex function $f:\R\to\R$,
  $x_1=-1$, $x_2=0$, $x_3=1$ and $l_1=l_3=1$, $l_2=2$, by
  Theorem~\ref{attaching} there exists a~polynomial $p\in\Pi_5$
  such that $p(-1)=f(-1)$, $p(0)=f(0)$, $p(1)=f(1)$ and
  $p(x)\le f(x)$ for $|x|>1$, $p(x)\ge f(x)$ for $0<|x|<1$. Observe
  that for a 5--convex function $f(x)=x^6$, $x\in\R$, this assertion
  is fulfilled by $p_1(x)=x^4$ and by $p_2(x)=x^2$.
\end{rem}
\begin{rem}
  The assumption $l_1=1$ if $x_1=\inf I$ is essential. For $n=k=1$ and
  $l_1=2$ Theorem~\ref{attaching} asserts that a~convex function
  $f:I\to\R$ has an affine support at a~point $x_1$. If $x_1$ is
  a~boundary point of $I$ it need not to be the case. Observe that
  a~convex function $f(x)=-\sqrt{1-x^2}$, $x\in[-1,1]$, has no affine
  support both at $x_1=-1$ and at $x_1=1$.
\end{rem}

\section{Some consequences of Theorem~\ref{attaching}}

We start with the support theorem for convex functions of an odd order.
\begin{cor}\label{support}
 Let $n\in\N$ be an odd number and $f:I\to\R$ be an $n$--convex
 function. Then for any $x_1\in\Int I$ there exists a~polynomial
 $p\in\Pi_n$ such that $p(x_1)=f(x_1)$ and $p\le f$ on~$I$.
\end{cor}
\begin{proof}
 Take in Theorem~\ref{attaching} $k=1$ and $l_1=n+1$.

\end{proof}
This result needs at least two comments.
\begin{enumerate}[1.]
 \item
  The support theorem for convex functions of an odd order was proved by
  Ger~\cite{G} with the additional assumption that the supported
  function is of the class~$\C^{n+1}$. However, Ger's result holds for
  functions defined on an open and convex subset of a~normed space.
  Notice at this place that if $f:I\to\R$ is $n$--convex then $f$ is
  of the class $\C^{n-1}$ on $\Int I$ (cf.~\cite{K,P}). Better
  regularity properties must be assumed (e.g. for $n=1$, $f(x)=|x|$
  is convex, continuous and not differentiable).
 \item
  The attaching method gives an answer to the question why convex
  functions of an even order need not to admit polynomial supports at
  every interior point of~$I$. Namely, for an $n$--convex function
  $f:I\to\R$ and $x_1\in\Int I$, if $n$ is an even number then the
  suitable interpolation polynomial $p$ (constructed as in the proof
  of Corollary~\ref{support}) fulfils by~\eqref{ineq:attaching} the
  inequality $p(x)\ge f(x)$, $x\in I$, $x<x_1$. It is easy to give an
  example: the 2--convex function $f(x)=x^3$ must not be supported on
  $\R$ by any quadratic polynomial (cf.~\cite{G}). However, there are
  some situations when convex functions of an even order do admit
  polynomial supports (cf. Corollary~\ref{cor_Rl} below).
\end{enumerate}
\begin{rem}
  Figures~\ref{Fig1} and~\ref{Fig2} show that the graph of a
  polynomial $p\in\Pi_n$ obtained by Theorem~\ref{attaching} may be
  situated on both sides of the graph of an $n$--convex function~$f$.
  Another possibility which may occur is that the graph of $p$ may
  be situated above the graph of $f$, contrary to the support
  property (see Corollaries~\ref{cor_L} and \ref{cor_Rr} below).
\end{rem}
\begin{cor}\label{cor_G}
 If $f:[a,b]\to\R$ is $(2n-1)$--convex and
 $x_1,\dots,x_n\in (a,b)$, then there exists a~polynomial
 $p\in\Pi_{2n-1}$ such that $p(x_i)=f(x_i)$, $i=1,\dots,n$,
 and $p\le f$ on $[a,b]$.
\end{cor}
\begin{proof}
 Assuming that $x_1<\dots<x_n$ use Theorem~\ref{attaching} for
 $2n-1$ instead of $n$, $k=n$ and $l_1=\dots=l_n=2$.
\end{proof}

\begin{cor}\label{cor_L}
 If $f:[a,b]\to\R$ is $(2n-1)$--convex and $x_1=a$,
 $x_2,\dots,x_n\in (a,b)$, $x_{n+1}=b$, then there exists a~polynomial
 $p\in\Pi_{2n-1}$ such that $p(x_i)=f(x_i)$, $i=1,\dots,n+1$, and
 $p\ge f$ on $[a,b]$.
\end{cor}
\begin{proof}
 Use Theorem~\ref{attaching} for $2n-1$ instead of $n$, $k=n+1$,
 $l_1=1$, $l_2=\dots=l_n=2$, $l_{n+1}=1$.
\end{proof}
\begin{cor}\label{cor_Rl}
 If $f:[a,b]\to\R$ is $2n$--convex, $x_1=a$,
 $x_2,\dots,x_{n+1}\in (a,b)$, then there exists a~polynomial
 $p\in\Pi_{2n}$ such that $p(x_i)=f(x_i)$, $i=1,\dots,n+1$,
 and $p\le f$ on $[a,b]$.
\end{cor}
\begin{proof}
 Use Theorem~\ref{attaching} for $2n$ instead of $n$, $k=n+1$ and
 $l_1=1$, $l_2=\dots=l_{n+1}=2$.
\end{proof}
\begin{cor}\label{cor_Rr}
 If $f:[a,b]\to\R$ is $2n$--convex, $x_1,\dots,x_n\in
 (a,b)$ and $x_{n+1}=b$, then there exists a~polynomial $p\in\Pi_{2n}$
 such that $p(x_i)=f(x_i)$, $i=1,\dots,n+1$, and $p\ge f$ on $[a,b]$.
\end{cor}
\begin{proof}
 Use Theorem~\ref{attaching} for $2n$ instead of $n$, $k=n+1$ and
 $l_1=\dots=l_n=2$, $l_{n+1}=1$.
\end{proof}

\section{Hadamard--type inequalities}

In this section we obtain some inequalities between the quadrature
operators and the integral approximated by them. The classical
inequality of this kind is the celebrated Hermite--Hadamard
inequality~\eqref{ineq:Hadamard}.

\subsection*{Orthogonal polynomials}

Let $w:[a,b]\to [0,\infty)$ be an integrable function such that
$\int_a^b w(x)dx>0$. The function~$w$ is called the~\emph{weight
function\/}. Then
\[
 \langle f,g\rangle_w:=\int_a^bf(x)g(x)w(x)dx
\] 
is the inner product in the space of all integrable functions
$f:[a,b]\to\R$. Performing for the sequence of monomials
$(1,x,x^2,\dots)$ the Gramm--Schmidt orthogonalization procedure we
obtain the sequence $(P_n)$ of polynomials orthogonal to each other on
$[a,b]$ with the weight~$w$ (i.e. with respect to the above inner
product). Let $P_n$ be the member of this sequence of degree~$n$. The
well--known results from numerical analysis (cf. e.g.~\cite{R,S})
state that the polynomial $P_n$ has $n$ distinct zeros belonging to
$(a,b)$.

\subsection*{Gauss quadratures}

Let $(P_n)$ be the sequence of polynomials orthogonal to each other on
$[a,b]$ with the weight function~$w$ and let $x_1,\dots,x_n$ be the
zeros of the polynomial~$P_n$. Furthermore, let
\begin{align*}
 w_i&:=\int_a^b\frac{P_n(x)w(x)}{(x-x_i)P'_n(x_i)}dx,
       \quad i=1,\dots,n,\\
 \G_n(f)&:=\sum_{i=1}^nw_if(x_i).
\end{align*}
It is well--known from numerical analysis
(cf. e.g.~\cite{BP2,R,S,W_G}) that the equation
\[
 \int_a^bf(x)w(x)dx=\G_n(f)
\]
holds for all polynomials belonging to $\Pi_{2n-1}$. If $[a,b]=[-1,1]$
and $w\equiv 1$ then $\G_n$ is the $n$--point Gauss--Legendre
quadrature (cf.~\cite{R,W_GL}).

\subsection*{Lobatto--type quadratures}

Let $(Q_n)$ be the sequence of polynomials orthogonal to each other on
$[a,b]$ with the weight function $(x-a)(b-x)w(x)$ and let
$x_1,\dots,x_{n-1}$ be the zeros of the polynomial~$Q_{n-1}$ (where
$Q_{n-1}$ is the member of this sequence of degree~$n-1$).
Furthermore, let
\begin{align*}
 w_0&:=\frac{1}{(b-a)Q_{n-1}^2(a)}\int_a^b Q_{n-1}^2(x)(b-x)w(x)dx,\\
 w_i&:=\frac{1}{(b-x_i)(x_i-a)}
 \int_a^b\frac{Q_{n-1}(x)(x-a)(b-x)w(x)}{(x-x_i)Q'_{n-1}(x_i)}dx,\\
 w_n&:=\frac{1}{(b-a)Q_{n-1}^2(b)}\int_a^b Q_{n-1}^2(x)(x-a)w(x)dx,\\
 \Lob_{n+1}(f)&:=w_0f(a)+\sum_{i=1}^{n-1}w_if(x_i)+w_nf(b).
\end{align*}
It is well--known from numerical analysis (cf. e.g.~\cite{BP2,R})
that the equation
\[
 \int_a^bf(x)w(x)dx=\Lob_{n+1}(f)
\]
holds for all polynomials belonging to $\Pi_{2n-1}$. If $[a,b]=[-1,1]$
and $w\equiv 1$ then $\Lob_{n+1}$ is the $(n+1)$--point Lobatto
quadrature (cf.~\cite{R,W_L}).

\subsection*{Inequalities for Gauss quadratures
             and Lobatto--type quadratures}

\begin{prop}\label{ineq:GL}
 If $f:[a,b]\to\R$ is $(2n-1)$--convex then
 \[
  \G_n(f)\le\int_a^b f(x)w(x)dx\le\Lob_{n+1}(f).
 \]
\end{prop}
\begin{proof}
 By Theorem~C $f$ is integrable on $[a,b]$. Let $x_1,\dots,x_n$ be the
 abscissas of the quadrature rule $\G_n$. By Corollary~\ref{cor_G} there
 exists a~polynomial $p\in\Pi_{2n-1}$ such that $p(x_i)=f(x_i)$,
 $i=1,\dots,n$, and $p\le f$ on $[a,b]$. Then $\G_n(p)=\G_n(f)$ and by
 $w\ge 0$
 \[
  \int_a^b p(x)w(x)dx\le\int_a^b f(x)w(x)dx.
 \]
 Since the quadrature $\G_n$ is precise for polynomials belonging to
 $\Pi_{2n-1}$, then
 \[
  \G_n(f)=\G_n(p)=\int_a^b p(x)w(x)dx\le\int_a^b f(x)w(x)dx.
 \]
 The second inequality we prove similarly taking as
 $x_1,\dots,x_{n+1}$ the abscissas of the quadrature rule $\Lob_{n+1}$
 and using Corollary~\ref{cor_L}.
\end{proof}

\subsection*{Radau--type quadratures}

Let $(P_n)$ be the sequence of polynomials orthogonal to each other on
$[a,b]$ with the weight function $(x-a)w(x)$ and let $x_1,\dots,x_n$
be the zeros of the polynomial~$P_n$. Furthermore, let
\begin{align*}
 w_0&:=\frac{1}{P_n^2(a)}\int_a^b P_n^2(x)w(x)dx,\\
 w_i&:=\frac{1}{(x_i-a)}
 \int_a^b\frac{P_n(x)(x-a)w(x)}{(x-x_i)P'_n(x_i)}dx,\;i=1,\dots,n,\\
 \Rad_{n+1}^l(f)&:=w_0f(a)+\sum_{i=1}^n w_if(x_i).
\end{align*}
If $[a,b]=[-1,1]$ and $w\equiv 1$ then $\Rad_{n+1}^l$ is the
$(n+1)$--point Radau quadrature (cf.~\cite{R,W_R}).
\par
Let $(Q_n)$ be the sequence of polynomials orthogonal to each other on
$[a,b]$ with the weight function $(b-x)w(x)$ and let $x_1,\dots,x_n$
be the zeros of the polynomial~$Q_n$. Furthermore, let
\begin{align*}
 w_i&:=\frac{1}{(b-x_i)}
 \int_a^b\frac{Q_n(x)(b-x)w(x)}{(x-x_i)Q'_n(x_i)}dx,\;i=1,\dots,n,\\
 w_{n+1}&:=\frac{1}{Q_n^2(b)}\int_a^b Q_n^2(x)w(x)dx,\\
 \Rad_{n+1}^r(f)&:=\sum_{i=1}^n w_if(x_i)+w_{n+1}f(b).
\end{align*}
It is well--known from numerical analysis
(cf. e.g.~\cite{BP2,R,W_R}) that the equation
\[
 \Rad_{n+1}^l(f)=\int_a^bf(x)w(x)dx=\Rad_{n+1}^r(f)
\]
holds for all polynomials belonging to $\Pi_{2n}$.

\subsection*{Inequalities for Radau--type quadratures}

\begin{prop}\label{ineq:R}
 If $f:[a,b]\to\R$ is $2n$--convex, then
 \[
  \Rad_{n+1}^l(f)\le\int_a^b f(x)w(x)dx\le\Rad_{n+1}^r(f).
 \]
\end{prop}
\begin{proof}
 The proof is similar to that of Proposition~\ref{ineq:GL}. For the
 first inequa\-lity use Corollary~\ref{cor_Rl} for the abscissas of the
 quadrature rule $\Rad_{n+1}^l$ and for the second one use
 Corollary~\ref{cor_Rr} for the abscissas of $\Rad_{n+1}^r$.

\end{proof}

\section{Comments}
\begin{enumerate}[1.]
 \item
  The inequalities of Propositions~\ref{ineq:GL} and~\ref{ineq:R} were
  earlier proved by Bessenyei and P\'ales \cite{BP1,BP2,BP_Cheb}.
 \item
  In \cite{BP2} these inequalities were proved for
  the weight function $w\equiv 1$ by the method of smoothing of convex
  functions of higher order. Namely, it is shown
  in~\cite[Theorem~5]{BP_Cheb} that for an $n$--convex
  function $f:I\to\R$ and for any compact subinterval $J\subset\Int I$
  there exists a~sequence of $n$--convex functions of the
  $\C^{\infty}$ class convergent uniformly to~$f$ on~$J$.
 \item
  In more recent paper~\cite{BP_Cheb} Hadamard--type inequalities for
  convex functions with respect to Chebyshev systems are given. The
  results are proved for any weight function. The method of the proof
  was based on integration of the determinant defining the convexity
  of this kind. The paper~\cite{BP1} contains the same results.
  However, some assumption present in~\cite{BP_Cheb} was removed. Both
  quoted papers do not contain any results of support--type.
 \item
  For some cases it is possible to give the inequalities of
  Hadamard--type which are better in some sense from the inequalities
  of Propositions~\ref{ineq:GL} and~\ref{ineq:R}. Some of them are
  presented in the next section.
\end{enumerate}

\section{Other Hadamard--type inequalities}

In this section we consider real functions defined on $[-1,1]$ and
the weight function $w\equiv 1$. In this setting
\begin{align*}
 \G_2(f)&=f\bigl(-\tfrac{\sqrt{3}}{3}\bigr)
         +f\bigl(\tfrac{\sqrt{3}}{3}\bigr),\\
 \Lob_4(f)&=\tfrac{1}{6}\bigl(f(-1)+f(1)\bigr)
           +\tfrac{5}{6}\Bigl(f\bigl(-\tfrac{\sqrt{5}}{5}\bigr)
           +f\bigl(\tfrac{\sqrt{5}}{5}\bigr)\Bigr).
\end{align*}
The abscissas of these quadrature rules are the zeros of suitable
orthogonal polynomials.

\subsection*{Remarks on even functions}

\begin{enumerate}[1.]
 \item
  If $f$ is an even function then
  $\int_{-1}^1 f(x)dx=2\int_0^1 f(x)dx$.
 \item
  If $f$ is an $n$--convex function and $n$ is an odd number then the
  function $f(-x)$ is also $n$--convex (cf.~\cite{P}). Then an even
  part of~$f$, i.e. the function $f_e(x)=\frac{f(x)+f(-x)}{2}$, is
  $n$--convex as well.
 \item
  Let $f$ be the integrable function. Then
  $\int_{-1}^1 f(x)dx=\int_{-1}^1 f_e(x)dx$. Indeed, since $f_e$ is
  an even function we have
  \begin{multline*}
   \int_{-1}^1 f_e(x)dx=2\int_0^1 f_e(x)dx
   =\int_0^1 \bigl(f(x)+f(-x)\bigr)dx\\
   =\int_0^1 f(x)dx+\int_0^1 f(-x)dx
   =\int_0^1 f(x)dx+\int_{-1}^0 f(t)dt=\int_{-1}^1 f(x)dx.
  \end{multline*}
 \item
  Fix $x_1,\dots,x_n\in (0,1]$ and for any function $f$
  define
  \[
   \T(f):=\alpha_0f(0)+\sum_{i=1}^n\alpha_i\bigl(f(x_i)+f(-x_i)\bigr).
  \]
  Then $\T(f)=\T(f_e)$. Namely,
  \[
   \T(f)=\alpha_0f(0)+\sum_{i=1}^n\alpha_i\cdot 2f_e(x_i)
   =\alpha_0f_e(0)+\sum_{i=1}^n\alpha_i\bigl(f_e(x_i)+f_e(-x_i)\bigr)
   =\T(f_e).
  \]
 \item
  Let $n$ be an odd positive integer. Because of the above remarks the
  inequalities of the form $\T_1(f)\le\int_{-1}^1 f(x)dx\le\T_2(f)$
  hold for any $n$--convex function~$f$ if and only if they hold for
  any $n$--convex and even function~$f$.
\end{enumerate}

\subsection*{Inequalities for Chebyshev quadrature}

Recall that the operator
\[
 \C(f):=\frac{2}{3}\Bigl(f\bigl(-\frac{\sqrt{2}}{2}\bigr)
       +f(0)+f\bigl(\frac{\sqrt{2}}{2}\bigr)\Bigr)
\]
is connected with the 3--point Chebyshev quadrature rule
(cf.~\cite{R,W_C}).
\begin{prop}\label{Cheb}
 If $f:[-1,1]\to\R$ is 3--convex then
 \[
  \G_2(f)\le\C(f)\le\int_{-1}^1 f(x)dx.
 \]
\end{prop}
\begin{proof}
 It is enough to prove the theorem for even functions.
 \begin{enumerate}[1.]
  \item
   By 3--convexity and~\eqref{ineq:nconv} $D(-v,-u,0,u,v;f)\ge 0$ for
   any $0<u<v\le 1$. Expanding this determinant by the last row we
   simply compute $v^2f(u)\le u^2f(v)+(v^2-u^2)f(0)$. For
   $u=\frac{\sqrt{3}}{3}$, $v=\frac{\sqrt{2}}{2}$ we obtain
   $\G_2(f)\le\C(f)$.
  \item
   By Theorem~\ref{attaching} (for $n=3$, $k=3$,
   $x_1=-\frac{\sqrt{2}}{2}$, $x_2=0$, $x_3=\frac{\sqrt{2}}{2}$,
   $l_1=l_2=1$, $l_3=2$) there exists a~polynomial $p\in\Pi_3$ such
   that
   $p\bigl(-\frac{\sqrt{2}}{2}\bigr)=f\bigl(-\frac{\sqrt{2}}{2}\bigr)$,
   $p(0)=f(0)$,
   $p\bigl(\frac{\sqrt{2}}{2}\bigr)=f\bigl(\frac{\sqrt{2}}{2}\bigr)$
   and $p\le f$ on~$[0,1]$. By Newton's Interpolation
   Formula~\eqref{eq:Newton}
   \begin{multline*}
    p(x)=f\bigl(\tfrac{\sqrt{2}}{2}\bigr)
        +\bigl[-\tfrac{\sqrt{2}}{2},0;f\bigr]
		 \bigl(x+\tfrac{\sqrt{2}}{2}\bigr)
        +\bigl[-\tfrac{\sqrt{2}}{2},0,\tfrac{\sqrt{2}}{2};f\bigr]
		 \bigl(x+\tfrac{\sqrt{2}}{2}\bigr)x\\
        +A\bigl(x+\tfrac{\sqrt{2}}{2}\bigr)x
		  \bigl(x-\tfrac{\sqrt{2}}{2}\bigr)
   \end{multline*}
  for some constant $A$. Computing these divided differences we can
  easily see that $2\int_0^1 p(x)dx=\C(f)$, whence
  $\C(f)\le2\int_0^1 f(x)dx=\int_{-1}^1 f(x)dx$.
 \end{enumerate}
\end{proof}

\subsection*{Inequalities for 5--convex functions}

Recall that the operator
$\Smp(f):=\frac{1}{3}\bigl(f(-1)+4f(0)+f(1)\bigr)$ is connected with
Simpson's quadrature rule (cf.~\cite{R,W_S}).
\begin{prop}\label{5conv}
 If $f:[-1,1]\to\R$ is 5--convex then
 \[
  \int_{-1}^1 f(x)dx
  \le\tfrac{2}{5}\Smp(f)+\tfrac{3}{5}\G_2(f)\le\Lob_4(f).
 \]
\end{prop}
\begin{proof}
 It is enough to prove the theorem for even functions.
 \begin{enumerate}[1.]
  \item
   By Theorem~\ref{attaching} ($n=5$, $k=5$, $x_1=-1$,
   $x_2=-\frac{\sqrt{3}}{3}$, $x_3=0$, $x_4=\frac{\sqrt{3}}{3}$,
   $x_5=1$, $l_1=l_2=l_3=1$, $l_4=2$, $l_5=1$) there exists
   a~polynomial $p\in\Pi_5$ such that $p(x_i)=f(x_i)$, $i=1,2,3,4,5$
   and $p\ge f$ on~$[0,1]$. Similarly as in the proof of
   Proposition~\ref{Cheb} we use Newton's Interpolation
   Formula~\eqref{eq:Newton} for the abscissas $x_1,x_2,x_3,x_4,x_5$
   and we compute
   $2\int_0^1p(x)dx=\frac{2}{5}\Smp(f)+\frac{3}{5}\G_2(f)$, from which
   the first inequality follows.
  \item
   To obtain the second inequality we also proceed similarly to the
   proof of Proposition~\ref{Cheb}. By 5--convexity
   and~\eqref{ineq:nconv}
   \[
    D\bigl(
	  -1,-\tfrac{\sqrt{3}}{3},-\tfrac{\sqrt{5}}{5},0,
	  \tfrac{\sqrt{5}}{5},\tfrac{\sqrt{3}}{3},1;f
	 \bigr)\ge 0.
   \]
   Expanding this determinant by the last row and performing some
   computations we get the desired inequality.
 \end{enumerate}
\end{proof}
Other inequalities between the quadrature operators can be found
in~\cite{W1}.

\section{Error bounds of quadrature rules}

Hadamard--type inequalities can be applied to estimate the errors
of quadrature rules. We illustrate this for the quadrature
$\T(f):=\frac{2}{5}\Smp(f)+\frac{3}{5}\G_2(f)$. Denote
$\I(f):=\int_{-1}^1 f(x)dx$.
\begin{prop}
 If $f\in\C^6\bigl([-1,1]\bigr)$ and
 $M:=\sup\Bigl\{\bigl|f^{(6)}(x)\bigr|\;:\;x\in [-1,1]\Bigr\}$, then
 $\bigl|\T(f)-\I(f)\bigr|\le\frac{M}{28350}$.
\end{prop}
\begin{proof}
 Let $g(x):=\frac{Mx^6}{6!}$. Then $g^{(6)}(x)=M$ and
 $\bigl|f^{(6)}(x)\bigr|\le g^{(6)}(x)$. Therefore $(g+f)^{(6)}\ge 0$
 and $(g-f)^{(6)}\ge 0$. By Theorem~D the functions $g+f$, $g-f$ are
 $5$--convex. By Proposition~\ref{5conv} we get
 \[
   \I(g+f)\le\T(g+f),\qquad\I(g-f)\le\T(g-f).
 \]
Since the operators $\T$ and $\I$ are linear then
 \[
   \I(g)-\T(g)\le\T(f)-\I(f),\qquad \T(f)-\I(f)\le\T(g)-\I(g).
 \]
 Hence $\bigl|\T(f)-\I(f)\bigr|\le\T(g)-\I(g)$. We conclude the proof
 by computing $\T(g)-\I(g)=\frac{M}{28350}$.
\end{proof}
The method presented above can be applied for other quadrature rules.
However, using it for Chebyshev, Gauss--Legendre, Lobatto, Radau and
Simpson's quadratures we obtain the error bounds known from numerical
analysis (cf.~\cite{R,W_C,W_GL,W_L,W_R,W_S}).

\bibliographystyle{amsplain}

\end{document}